\documentclass[final]{siamltex}
\usepackage{amsmath}
\usepackage{amsfonts}
\usepackage{color}
\usepackage[latin1]{inputenc}
\usepackage{algpseudocode}
\usepackage{algorithm}
\usepackage{caption}
\usepackage{mathrsfs}
\usepackage{subcaption}
\usepackage{lscape}
\usepackage{mwe}
\usepackage{enumerate}

\def\bxi{{\boldsymbol{\xi}}}
\def\bth{{\boldsymbol{\theta}}}
\def\bu{{\bf{u}}}
\usepackage[update,prepend]{epstopdf} 
\usepackage{graphicx}

\newcommand{\be}{\begin{eqnarray}}
\newcommand{\e}{\end{eqnarray}}
\newcommand{\bes}{\begin{eqnarray*}}
\newcommand{\es}{\end{eqnarray*}}
\newcommand{\beq}{\begin{equation}}
\newcommand{\eeq}{\end{equation}}

\newcommand{\R}{\ensuremath{\mathbb{R}}}

\usepackage[textsize=small]{todonotes}
\title{Novel Deep neural networks for solving \\ Bayesian statistical inverse 
problems \thanks{HA, AO, and DV are partially supported by NSF grants DMS-1818772, DMS-1913004, 
the Air Force Office of Scientific Research (AFOSR) under Award NO: FA9550-19-1-0036, and 
Department of Navy, Naval PostGraduate School under Award NO: N00244-20-1-0005. HCE is partially 
supported by NSF grant DMS-1819115.}}

\begin{document}

\author{Harbir Antil\thanks{Department of Mathematical Sciences
			   and The Center for Mathematics and Artificial Intelligence, 
                            George Mason University,
                            Fairfax, VA 22030, USA.
                            \texttt{hantil@gmu.edu}.
			   }                           
\and
Howard C. Elman\thanks{Department of Computer Science and Institute for 
                     Advanced Computer Studies, University of Maryland, 
                     College Park, MD 20742, USA.
                     \texttt{helman@umd.edu}.
		    }
\and 
Akwum  Onwunta\thanks{Department of Mathematical Sciences
			   and The Center for Mathematics and Artificial Intelligence, 
                            George Mason University,
                            Fairfax, VA 22030, USA. 
                          \texttt{aonwunta@gmu.edu}.                          
			 }     
\and			                 
Deepanshu Verma\thanks{Department of Mathematical Sciences
			   and The Center for Mathematics and Artificial Intelligence, 
                            George Mason University,
                            Fairfax, VA 22030, USA. 
                            \texttt{dverma2@gmu.edu}.                          
			 }                    
}

\maketitle
\begin{abstract}
We consider the simulation of Bayesian statistical inverse problems 
governed by large-scale linear and nonlinear partial differential
equations (PDEs).  Markov chain Monte Carlo (MCMC) algorithms 
are standard techniques to solve such problems. However, MCMC
techniques are computationally challenging as they require several
thousands of forward PDE solves. The goal of this paper is to introduce 
a fractional deep neural network based approach for the forward solves within 
an MCMC routine. Moreover, we discuss some approximation error estimates and  
illustrate the efficiency of our approach via several numerical examples.
\end{abstract}



\section{Introduction}
\label{sec1}

Large-scale  statistical inverse  problems governed by partial differential equations (PDEs) are increasingly found in different areas of computational science and engineering \cite{BSHL2014, Bui14, CS07, Flath11, KS05}.  
The basic use of this class of problems  is to recover certain physical quantities from limited and noisy observations. They are generally computationally demanding and can be solved using the Bayesian framework \cite{Bui14,  EO19, Flath11}. In the  Bayesian
approach to statistical inverse problems, one models  the solution as a
posterior distribution of the unknown parameters conditioned on the observations. Such mathematical formulations are often ill-posed and  regularization is introduced in the model via prior information. The
posterior distribution completely characterizes the uncertainty in the model. Once it has been computed, statistical quantities of interest can then be obtained from the posterior distribution. Unfortunately, many problems of practical relevance do not admit analytic representation of the posterior distribution. Thus, the posterior distributions are generally sampled using  Markov chain Monte Carlo (MCMC)-type schemes.

We point out here that a large number of samples are generally
needed by MCMC methods to obtain convergence. This makes MCMC methods computationally intensive to simulate Bayesian inverse  problems. Moreover, 
for  statistical inverse problems governed by  PDEs, one needs to solve the forward problem corresponding to each MCMC sample. This task further increases the computational complexity of the problem, especially if the governing PDE is a large-scale nonlinear PDE \cite{EO19}.
Hence, it is reasonable to  construct a computationally cheap surrogate 
 to replace the forward model solver \cite{LM14}. Problems that involve  large numbers of input
parameters often lead to  the so-called {\it curse of dimensionality}.  Projection-based reduced-order models  e.g., reduced basis methods and the discrete empirical interpolation method (DEIM) are  typical examples of dimensionality
reduction methods for tackling parametrized PDEs \cite{HAntil_MHeinkenschloss_DCSorensen_2013a, BMNP04, CS10, EF17, Forstall2015,JSHesthaven_GRozza_BStamm_2016a,AQuarteroni_AManzoni_FNegri_2016a}. However, they are intrusive by design in the sense that they do not allow reuse of existing codes in the forward solves, especially for nonlinear forward models. To mitigate this computational issue, the goal of this work is to demonstrate the use of deep neural networks (DNN) to construct surrogate models  specifically for nonlinear parametrized PDEs governing statistical Bayesian inverse  problems. In particular, we will focus on  fractional Deep Neural Networks (fDNN) which have been recently introduced   and  applied to classification problems \cite{Antil0420}.

The use of  DNNs for surrogate modeling in the framework of   PDEs has received increasing attention in recent years, see e.g., \cite{Chen:20, GRPK19, HJE18, ZLi20, PLK19,  RPK19, TB18, YANG2021109913, ZZ18}  and the references therein. Some of these references also cover the type of Bayesian inverse problems considered in our paper. To accelerate or replace computationally-expensive PDE solves, a DNN is trained to approximate the mapping from parameters in the PDE to observations obtained from its solution. In a supervised learning setting, training data consisting of inputs and outputs are available and the learning problem aims 
at tuning the weights of the DNN. To obtain an effective surrogate model, it is desirable to train the DNN to a high accuracy. 

We consider a different approach than the aforementioned works. We propose novel dynamical systems based neural 
networks which allow connectivity among all the network layers.  
Specifically, we consider a {fractional} DNN technique  recently proposed in \cite{Antil0420} in the framework of classification problems. We note that 
\cite{Antil0420} is motivated by \cite{Lars2020}. Both papers are in the spirit of the push to develop rigorous mathematical models for the analysis and understanding of  DNNs. The idea is to consider DNNs as dynamical systems. More precisely,
in \cite{Benning19, Lars2020, HR17}, DNN is thought of as an optimization problem constrained by a discrete ordinary differential equation (ODE). 
As pointed out in \cite{Antil0420},  designing the DNN solution algorithms
at the continuous level has the appealing advantage of  architecture independence; in other words, the number of optimization iterations remains the same even if the number of layers is increased. Moreover, \cite{Antil0420} specifically considers continuous fractional ODE constraint. Unlike standard DNNs, the resulting fractional DNN allows the network to access historic information of input and gradients across all subsequent layers since all the layers are connected to one another.

We demonstrate that our fractional DNN  leads to a significant reduction in the overall 
computational time used by MCMC algorithms to solve inverse problems, while maintaining accurate evaluation of the relevant statistical quantities. 
We note here that, although the papers \cite{ZLi20, TB18, ZZ18}  consider DNNs for inverse problems governed by PDEs, the DNNs used in these studies do not, in general, have the optimization-based formulation  considered here. 

The remainder of the paper is organized as follows: In section~\ref{probform} we state the generic parametrized PDEs under consideration, as well as  discuss the well-known surrogate modeling approaches. Our main DNN algorithm to approximate parame\-trized PDEs is provided in section~\ref{dnn}.  We first discuss a ResNet architecture to approximate parametrized PDEs, which is followed by our fractional ResNet (or fractional DNN) approach to carry out the same task. A brief discussion on error estimates has been provided. Next, we discuss the application of fractional DNN to statistical Bayesian inverse problems in section \ref{inverse_problems}. We conclude the paper with several illustrative numerical examples in section~\ref{Numex}.

\section{Parametrized PDEs}
\label{probform}
In this section, we present an abstract formulation of the forward problem as a discrete  partial differential equation (PDE) depending on parameters. 
We are interested in the following   model which represents (finite element, finite difference or finite volume) discretization of a (possibly nonlinear) parametrized PDE
\be
\label{disc1}
F(\bu(\bxi);\bxi) = 0, 
\e
where $\bu(\bxi)\in \mathcal{U}\subset \mathbb{R}^{N_x}$  and  
 $\bxi \in \mathcal{P}\subset\mathbb{R}^{N_{\bxi}}$ denote  the solution of the PDE  and the  parameter in the model, respectively. Here, $\mathcal{P}$ denotes the parameter domain and $\mathcal{U}$ the solution manifold. For a fixed parameter $\bxi \in \mathcal{P}$, we seek  the solution $\bu(\bxi)\in \mathcal{U}$. In other words, we have the functional relation given by the parameter-to-solution map 
\begin{equation}
	\label{eq:Phi}
 	\bxi\mapsto\Phi(\bxi)\equiv \bu(\bxi).
\end{equation}
 
Several  processes in computational sciences and engineering are modelled via parameter-dependent 
 PDEs, which when discretized are of the form (\ref{disc1}), for instance, in  viscous flows governed by Navier-Stokes equations, which is parametrized by the Reynolds number \cite{ESW14}. Similarly,  in  unsteady natural
convection problems  modelled via Boussinesq equations,  the Grashof or Prandtl numbers are  important parameters  \cite{AQuarteroni_AManzoni_FNegri_2016a}, etc.
Approximating the high-fidelity solution of (\ref{disc1}) can be done with  nonlinear solvers, such as Newton or Picard methods combined with Krylov solvers \cite{ESW14}. However, computing solutions to (\ref{disc1}) can become prohibitive, especially when they are required for many parameter values and $N_{\bxi}$ is  large. Besides, a relatively large $N_x$ (resulting from a really fine mesh in the discretization of the PDE) yields large (nonlinear) algebraic systems  which are computationally expensive to solve and may also lead to huge storage requirements. This is, for instance, the case in Bayesian inference problems governed by PDEs where several forward solves are required to adequately sample posterior distributions through  MCMC-type schemes \cite{BGLFS17, Martin12}.
As a result of the above computational challenges, it is reasonable to replace the high-fidelity model by a  surrogate model which is relatively easy to evaluate.  

In what follows, we  discuss two classes of surrogate models, namely: reduced-order models and deep learning models. 

\subsection{Surrogate Modeling}
\label{surrmod}
Surrogate  models  are cheap-to-evaluate models  designed to replace  computationally costly models. The major advantage  of  surrogate models is that approximate  solution of the models can be easily evaluated at any new parameter instance with  minimal loss of accuracy,  at  a cost independent of the dimension of the  original high-fidelity problem.

A popular class of surrogate models are the reduced-order models (ROMs), see e.g., \cite{HAntil_MHeinkenschloss_DCSorensen_2013a, BMNP04, Cha2011, EF15, EF17, Forstall2015,JSHesthaven_GRozza_BStamm_2016a,AQuarteroni_AManzoni_FNegri_2016a}. Typical examples of ROMs include the reduced basis method \cite{AQuarteroni_AManzoni_FNegri_2016a}, proper orthogonal decomposition \cite{Cha2011, GMNP07, JSHesthaven_GRozza_BStamm_2016a, AQuarteroni_AManzoni_FNegri_2016a} and the discrete empirical
interpolation method  (DEIM) and its variants \cite{HAntil_MHeinkenschloss_DCSorensen_2013a, BMNP04, Cha2011, EF17, Forstall2015}. A key feature of the ROMs is that they use the so-called {\it offline-online paradigm}. 
The offline step essentially constructs the low-dimensional approximation to the solution space; this approximation is generally known as the {\it reduced basis}. The online step uses the reduced basis to solve a smaller reduced problem. The resulting reduced solution  accurately approximates the solution of the original problem.

Deep neural network (DNN) models constitute another class of surrogate models which are well-known for their high approximation capabilities. The basic idea of DNNs is to approximate  multivariate functions through a set of layers of increasing complexity \cite{GBC2016}. Examples of DNNs for surrogate modeling include Residual Neural Network (ResNet) \cite{HR17, HXRS16}, physics-informed neural network (PINNs) \cite{RPK19} and fractional DNN \cite{Antil0420}.
 
Note that the ROMs require system solves and they are highly intrusive especially for nonlinear problems \cite{HAntil_MHeinkenschloss_DCSorensen_2013a, Cha2011, EF17, EO19}. In contrast, the DNN approach is fully non-intrusive, which is essential for legacy codes. Although rigorous error estimates for ROMs under various assumptions have been well studied \cite{AQuarteroni_AManzoni_FNegri_2016a}, 
 we like the advantage of DNN  being nonintrusive, but recognize that error analysis is not yet as strong.

 In this work, we propose a surrogate model based on  the combination of POD and fractional DNN. Before we discuss our proposed model, we  first  review POD.
 
\subsection{Proper orthogonal decomposition}
\label{pod}
For sufficiently large   $N_s\in \mathbb{N}$,  suppose that   $\mathbb{E}:=\{\bxi_1,\bxi_2,\cdots,\bxi_{N_s}\}$ is a set of parameter samples with $\bxi_i\in\mathbb{R}^{N_{\bxi}}$, and  $\{\bu(\bxi_1),\bu(\bxi_2),\cdots,\bu(\bxi_{N_s})\}$ the corresponding snapshots (solutions of the model (\ref{disc1}), with $\bu_i\in\mathbb{R}^{N_x}$). Here, we assume that $\mbox{span}\{\bu(\bxi_1),\bu(\bxi_2),\cdots,\bu(\bxi_{N_s})\}$  sufficiently approximates the space of all possible  solutions of (\ref{disc1}).
Next, we denote by 
\be
\mathbb{S}=[\bu_1|\bu_2|\cdots| \bu_{N_s}]\in \R^{N_{x}\times N_s}
\e
the matrix whose columns are the solution snapshots. Then, the singular value decomposition (SVD) of $\mathbb{S}$ is given by 
\be
\label{svd}
\mathbb{S} = \widetilde{\mathbb{V}}\Sigma \mathbb{W}^T,
\e
where $\widetilde{\mathbb{V}}\in\mathbb{R}^{N_x\times r}$ and $\mathbb{W}\in\mathbb{R}^{N_s\times r}$  are orthogonal matrices called the left and right singular vectors, and $r\leq \min\{N_s,N_x\}$ is the rank of $\mathbb{S}$. Thus, $\widetilde{\mathbb{V}}^T \widetilde{\mathbb{V}}=I_r,\;\; \mathbb{W}^T\mathbb{W}=I_r,$  and $\Sigma=\mbox{diag}(\rho_1,\rho_2,\cdots,\rho_r)\in\mathbb{R}^{r\times r},$ where  $\rho_1\geq \rho_2\geq  \cdots \geq \rho_r \geq 0$ are the singular values of $\mathbb{S}.$ 

Now, denote by $\mathbb{V} \subset \widetilde{\mathbb{V}}$ the first $k\leq r$ left singular vectors of   $\mathbb{S}.$ Then, the columns of $\mathbb{V}\in \mathbb{R}^{N_x\times k} $ form a POD basis of dimension $k.$ 
According to the Schmidt-Eckart-Young theorem \cite{EY36, GV96}, the POD basis $\mathbb{V} $ 
 minimizes, over all
possible $k$-dimensional orthonormal bases $\mathbb{Z} \in \mathbb{R}^{N_x\times k}$, the sum of the squares of the errors between each snapshot vector $\bu_i$ and its projection onto the subspace spanned by $\mathbb{Z}$. More precisely,
\be
\label{EYthm}
\sum_{i=1}^{N_s}||\bu_i - \mathbb{V} \mathbb{V}^T\bu_i ||_2^2 = \min_{\mathbb{Z}\in\mathcal{Z}}\sum_{i=1}^{N_s}||\bu_i - \mathbb{Z} \mathbb{Z}^T\bu_i ||_2^2 = \sum_{i=k+1}^{r} \rho_i,
\e
where $\mathcal{Z}:=\{\mathbb{Z} \in \mathbb{R}^{N_x\times k}: \mathbb{Z}^T\mathbb{Z} =I_{k}\}.$
Note from (\ref{EYthm}) that the POD basis $\mathbb{V} $ solves a 
least squares minimization problem, which guarantees that the approximation error is controlled by the 
singular values.

For every $\bxi$, we then approximate the continuous solution $\bu(\bxi)$ as $\bu(\bxi) \approx \mathbb{V}\widehat{\bu}(\bxi)$, where 
$\widehat{\bu}(\bxi)$ solves the reduced problem 
\be
\label{disc1_rom}
	\mathbb{V}^T F(\mathbb{V}\widehat{\bu}(\bxi);\bxi) = 0 .
\e 
Notice that, in some reduced-modeling techniques such as  DEIM, additional steps are needed to fully reduce the dimensionality of the problem \eqref{disc1_rom}. Nevertheless, one still needs to solve a nonlinear (reduced) system like \eqref{disc1_rom} to evaluate $\widehat{\bu}$. 

We conclude this section by emphasizing that the above approach is ``linear" because $\bu(\bxi) \approx \widehat{\Phi}(\bxi)$, where $\widehat{\Phi} (\bxi)= \mathbb{V}\widehat{\bu}(\bxi)$ is an approximation of the map $\Phi(\bxi)$ given in \eqref{eq:Phi}.

\section{Deep neural network}
\label{dnn} 

The  DNN approach to modeling  surrogates  produces a nonlinear approximation $\widehat{\Phi}$ of the input-output map 
$\Phi:\mathbb{R}^{N_{\bxi}}  \rightarrow \mathbb{R}^{N_x}$ given in \eqref{eq:Phi}, where  $\widehat{\Phi}$  depends implicitly on  a set of parameters $\bth  \in   \mathbb{R}^{N_{\theta}}$ configured as a layered set of latent variables that must be trained. We represent this dependence using the notation $\widehat{\Phi}(\bxi;{\bth}).$ 
In the context of PDE surrogate modeling, training  a DNN requires 
a data set $(\mathbb{E},\mathbb{S}),$ where  the parameter samples $\bxi_j\in\mathbb{E}$ are the inputs   and the corresponding snapshots ${\bu_j}\in \mathbb{S}$ are  targets; training then consists of constructing 
$\bth$ so that the DNN $\widehat{\Phi}(\bxi_j;\bth)$ matches $\bu_j$ for each $\bxi_j.$ 
This matching is determined by a \emph{loss functional}.
The learning problem therefore involves computing the optimal parameter  $\bth$ that minimizes the loss functional and satisfies 
 $\widehat{\Phi}({\bxi}_j; \bth)\approx {\bu_j}.$ The ultimate goal is that this  approximation also holds with the same optimal parameter $\bth$  for a different data set; in other words, for $\bxi\notin \mathbb{E}, $ we take $\widehat{\Phi}(\bxi; \bth)$ to represent a good approximation to $\Phi(\bxi).$

 Deep learning can be either supervised or unsupervised depending on the data set used in training.  In the supervised learning technique for DNN, all the input samples ${\bxi_j}$  are  available for all the corresponding samples of the targets ${\bu_j}.$  In contrast, the unsupervised learning framework does not require all the outputs to accomplish the training phase. 
 In this paper, we adopt the supervised learning approach  to model our surrogate and apply it to  Bayesian inverse problems.  In particular, we shall discuss Residual neural network (ResNet)\cite{HR17, HXRS16} and Fractional DNN   \cite{Antil0420} in the context of PDE surrogate modeling.
   
\subsection{Residual neural network}
\label{subsec2a}
\begin{figure}[htb]
	\centering
		\includegraphics[width=0.48\textwidth]{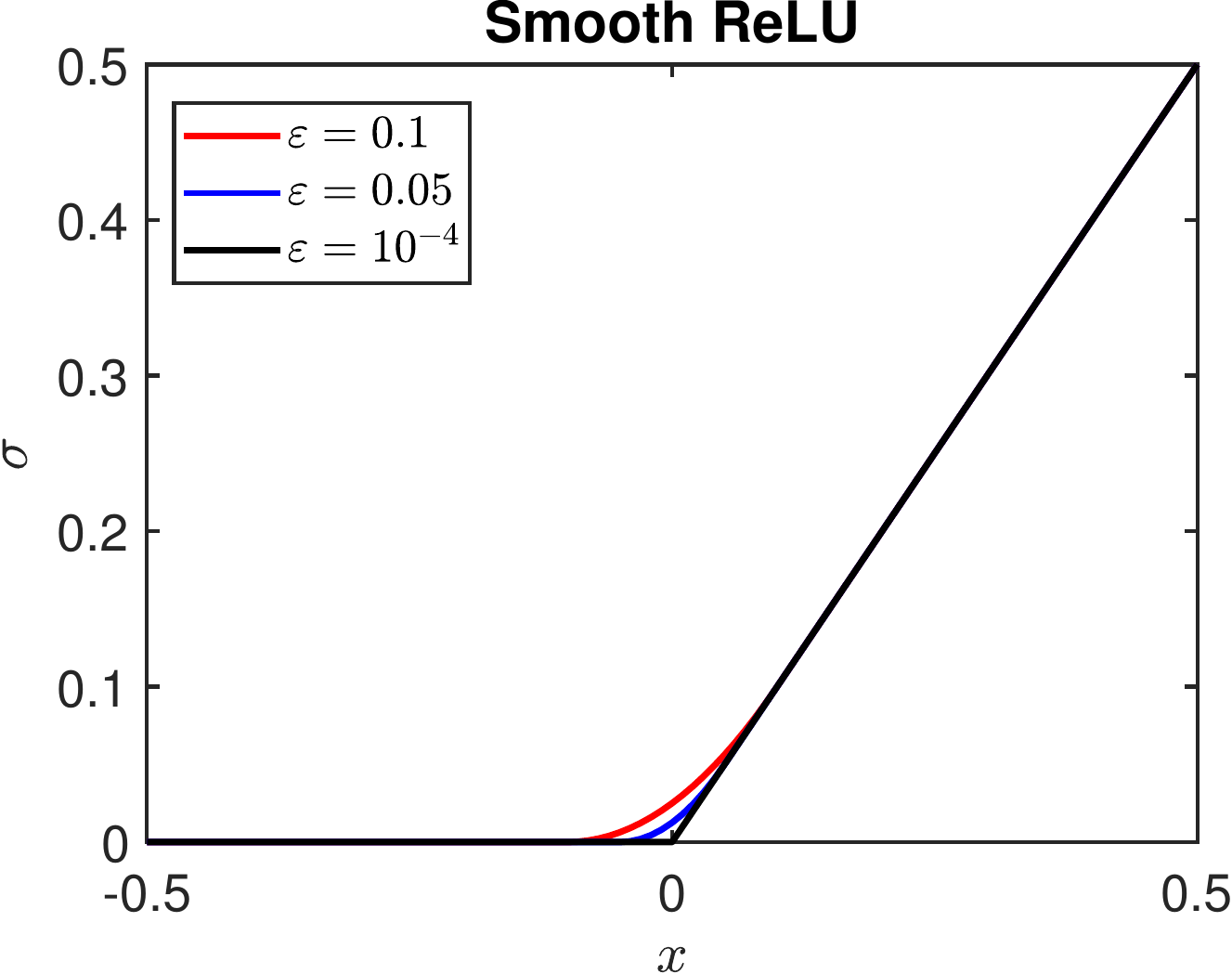}
		\includegraphics[width=0.48\textwidth]{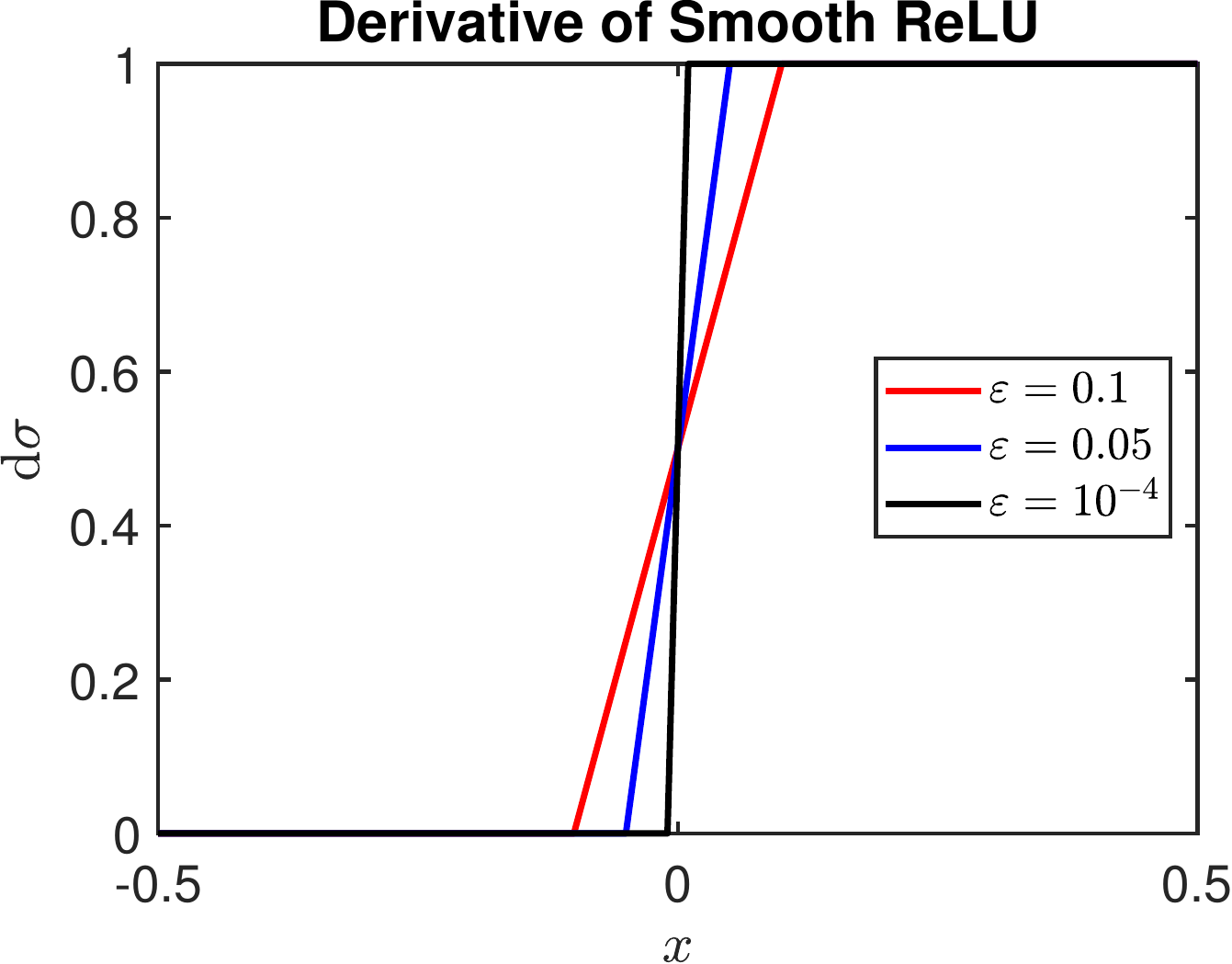}
		\caption{\label{fig:smooth_relu}Smooth ReLU and  derivative for different values of $\varepsilon$.}
\end{figure}

The residual neural network (ResNet) model was originally proposed in \cite{HXRS16}. For a given input datum $\bxi$, ResNet approximates $\Phi$ 
through the following recursive expression 
\begin{align}\label{eq:RNN}
\begin{aligned}
\phi_1 &= \sigma(W_0{\bxi} + {\bf b}_0), \\
\phi_{j} &= \phi_{j-1} + h \sigma(W_{j-1}\phi_{j-1} + {\bf b}_{j-1}), \quad 2\leq j\leq L-1 , \\
\phi_L   &= W_{L-1} \phi_{L-1} , 
\end{aligned}
\end{align}
where $\{W_j,{\bf b}_j\}$ are the weights and the biases, $h>0$ is the stepsize and $L$ is the number of layers and
 $\sigma$ is an activation function which is applied element-wise on its arguments. Typical examples of   $\sigma$ include the hyperbolic tangent function,   the logistic function or the rectified linear unit (or ReLU) function. ReLU is a nonlinear function given by $\sigma(x)=\max\{x,0\}.$ In this work, we use a smoothed ReLU function defined, for $\varepsilon>0,$  as 
\be
\label{smrelu}
	\sigma(x):=\begin{cases}
		x,& x>\varepsilon\\
		0,& x<-\varepsilon\\
		\frac{1}{4 \varepsilon}x^2 +\frac{1}{2}x+\frac{\varepsilon}{4} ,& -\varepsilon \le x\le \varepsilon.
	\end{cases} 
\e
Note that as $\varepsilon \to 0$, smooth ReLU approaches ReLU, see Figure~\ref{fig:smooth_relu}.

It follows from (\ref{eq:RNN}) that 
	 \begin{equation}
		{\widehat{\Phi}}(\bxi;\bth) = \phi_{L}(\bxi)=W_{L-1}\Big( 
		\big(I +h (\sigma \circ \mathcal{K}_{{L-1}})\big) \circ \dots \circ \big(I +h (\sigma \circ \mathcal{K}_{{1}})\big) \circ  \mathcal{K}_{0}
		\Big)(\bxi),
	\end{equation}
	where $\mathcal{K}_{j}({\bf y})=W_{j}\,{\bf y} + {\bf b}_{j},$ for all $j=0,\hdots,L-1$ and for any $\bf y$.
		
To this end, the following two critical questions naturally come to mind:
	\begin{enumerate}[(a)]
		\item How well does $\widehat{\Phi}(\bxi;\bth)$ approximate $\Phi(\bxi)$?
		\item How do we determine $\bth$?
	\end{enumerate}

	We will address (b) first and leave (a) for a discussion provided in section~\ref{erroranaly}.
Now, setting $\bth =\{W_j,{\bf b}_j\},$  it follows that the problem of
approximating $\Phi$ via ResNet is essentially the problem of learning the unknown parameter $\bth.$ More specifically, the learning problem is the solution of the minimization problem \cite{HR17}: 
\be
\label{optprb}
\min_{\bth} \mathcal{J(\bth ; \bxi,\bu)}
\e
subject to constraints \eqref{eq:RNN}, where $\mathcal{J(\bth, \bxi,\bu)}$ is suitable loss function. 

In this work, we consider the mean squared error, together with a regularization term, as our loss functional:
\be
\label{lossf}
 \mathcal{J}(\bth;\bxi, \bu)= \frac{1}{2N_s}\sum_{j=1}^{N_s}\| \hat{\Phi}(\bxi_j;\bth)-\bu_j \|^2_2 + \frac{\lambda}{2}||\bth||_2^2,
\e
where $\lambda$ is the regularization parameter. 
Due to its highly non-convex nature, this is a very difficult optimization problem. Indeed, a search over a high dimensional parameter space  for the global minimum of a non-convex function can be  intractable. The current state-of-the-art approaches to solve these optimization problems are based on the stochastic gradient descent method \cite{ZLi20, TB18, ZZ18}.

As pointed out in, for instance, \cite{HR17}, the middle equation in expression (\ref{eq:RNN}) mimics the forward Euler discretization of a nonlinear differential equation: 
\bes
\label{diffeq}
d_t\phi(t) &=& \sigma(W(t)\phi(t)  + b(t)), \;\; t\in (0,T] ,\\
\phi(0) &=& \phi_0.
\es

It is known  that standard DNNs are prone to vanishing gradient problem \cite{veit2016residual}, leading to loss of information during the training phase. ResNets do help, but more helpful is the so-called DenseNet \cite{huang2017denseNet}. Notice that in a typical DenseNet, each layer takes into 
account all the previous layers. However, this is an ad hoc approach with no rigorous
mathematical framework. Recently in \cite{Antil0420}, a mathematically rigorous approach
based on fractional derivatives has been introduced. This ResNet is called \emph{Fractional DNN.}
In \cite{Antil0420},
the authors numerically establish that fractional derivative based ResNet outperforms  ResNet in overcoming the vanishing gradient 
problem. This is not surprising because fractional derivatives allow connectivity among all the layers.  
Building rigorously on the idea of \cite{Antil0420}, we proceed to present the fractional DNN surrogate model for the discrete PDE in \eqref{disc1}.

\subsection{Fractional deep neural network}
\label{subsec2b}
As  pointed out in the previous section, the fractional DNN approach is based on the replacement of the standard ODE constraint in the learning problem by a fractional differential equation. To this end, we first introduce the definitions and concepts on which we shall rely to discuss fractional DNN. 

Let $u:[0,T]\rightarrow \mathbb{R}$  be an absolutely continuous function and assume $\gamma\in (0,1)$. Next, consider the following fractional differential equations 
\begin{equation} \label{eq:eq1}
d_t^{\gamma}u(t)=f(u(t)), \quad u(0)=u_0,
\end{equation}
and
\begin{equation} \label{eq:eq2}
d_{T-t}^{\gamma}u(t)=f(u(t)), \quad u(T)=u_T. 
\end{equation}
Here, $d_t^{\gamma}$ and  $d_{T-t}^{\gamma}$  denote left and right Caputo derivatives, respectively \cite{Antil0420}.

Then, setting $u(t_j)=u_j, $ and using the $L^1$-scheme (see e.g., \cite{Antil0420}) for the discretization of (\ref{eq:eq1}) and (\ref{eq:eq2}) yields, respectively,
\begin{alignat}{3} \label{eq:disc_cCaputo}
u_{j+1}= u_j - \sum_{k=0}^{j-1}   a_{j-k}\: \left(u_{k+1} - u_k \right)  + h^{\gamma}\Gamma(2-\gamma)f(u_j) ,
\quad  j = 0,...,L-1, \;
\end{alignat}
and
\begin{alignat}{3} \label{eq:disc_cRCaputo}
u_{j-1}= u_j + \sum_{k=j}^{L-1}   a_{k-j}\: \left(u_{k+1} - u_k \right)  - h^{\gamma}\Gamma(2-\gamma)f(u_j) , 
\quad  j = L,...,1,
\end{alignat}
where $h>0$ is the step size,  $\Gamma(\cdot)$ is the Euler-Gamma function, and  
\begin{equation}\label{a_k}
a_{j-k}=(j+1-k)^{1-\gamma}-(j-k)^{1-\gamma}.
\end{equation}

After this brief overview of the fractional derivatives, we are ready to introduce our fractional DNN (cf.~\eqref{eq:RNN}) 
\begin{alignat}{5}\label{eq:fRNN}
\phi_1 &= \sigma(W_0 \phi_0 + {\bf b}_0); \quad \phi_0=\bxi, \nonumber \\
\phi_{j} &= \phi_{j-1} - \sum_{k=0}^{j-2}   a_{j-1-k}\: \left(\phi_{k+1} - \phi_k \right)  + h^{\gamma}\Gamma(2-\gamma)&\sigma(W_{j-1}\phi_{j-1} + {\bf b}_{j-1}),\\& &  2\leq j\leq L-1 ,\nonumber \\
\phi_L   &= W_{L-1} \phi_{L-1} . \nonumber
\end{alignat}
Our learning problem then amounts to 
\be
\label{eq:fDNNJ}
\min_{\bth} \mathcal{J}(\bth ; \bxi,\bu)
\e
subject to constraints \eqref{eq:fRNN}. Notice that the middle equation in \eqref{eq:fRNN} mimics the $L^1$-in time discretization of 
the following nonlinear fractional differential equation
\begin{equation*}
\begin{aligned}
d_t^\gamma \phi(t) &= \sigma(W(t)\phi(t)  + b(t)), \;\; t\in (0,T] ,\\
\phi(0) &= \phi_0.
\end{aligned}
\end{equation*}

There are two ways to approach the constrained optimization problem \eqref{eq:fDNNJ}. The first approach is the so-called reduced approach, where we eliminate the constraints \eqref{eq:fRNN} and consider the minimization problem \eqref{eq:fDNNJ} only in terms of 
$\bth$.  The resulting problem can be solved by using a gradient based method such as BFGS, see e.g., \cite[Chapter 4]{CTK99}. During every step of the gradient method, one needs to solve the state equation \eqref{eq:fRNN} and an adjoint equation. These two solves enable us to derive an expression of the gradient of the reduced loss function with respect to $\bth$. Alternatively, one can derive the gradient and adjoint equations by using the Lagrangian approach. It is well-known that the gradient with respect to $\bth$ for both approaches coincides, see \cite[pg.~14]{MR3839730} for instance. We next illustrate  how to evaluate this gradient using the Lagrangian approach. 
	
The Lagrangian functional associated with the discrete constrained optimization problem \eqref{eq:fDNNJ} is given by
	\begin{multline}\label{eq:lagfdnn}
		\mathcal{L}(\bu,\bth,\boldsymbol{\psi}) := \mathcal{J(\bth ; \bxi,\bu)} +	\langle \phi_1 -\sigma(W_0 \phi_0 + {\bf b}_0), \psi_1  \rangle \\
		 + \sum_{j=2}^{L-1} \Bigg\langle \phi_{j} - \phi_{j-1} + \sum_{k=0}^{j-2}   a_{j-1-k}\: \left(\phi_{k+1} - \phi_k \right) +  
		  h^{\gamma}\Gamma(2-\gamma)\sigma(W_{j-1}\phi_{j-1} + {\bf b}_{j-1}), \psi_j  \Bigg\rangle \\
		 + \langle \phi_L - W_{L-1} \phi_{L-1}, \psi_L  \rangle	,
	\end{multline}
	where $\psi_j$'s are the Lagrange multi1pliers, also called adjoint variables,  corresponding to \eqref{eq:fRNN} and $\langle \cdot, \cdot \rangle$ is the inner product on $\mathbb{R}^{N_x}$.
	
Next we write the state and adjoint equations fulfilled at a stationary point of the Lagrangian $\mathcal{L}$. In addition, we state the derivative of $\mathcal{L}$ with respect to the design variable $\bth$. 
	\begin{subequations}\label{eq:oc_fDNN}
		\begin{enumerate}[(i)]
			\item State Equation. 
			\begin{align} \label{eq:oc_state}	
			\begin{aligned}
				\phi_1 &= \sigma(W_0{\bxi} + {\bf b}_0),  \\
				\phi_{j} &= \phi_{j-1} - \sum_{k=0}^{j-2}   a_{j-1-k}\: \left(\phi_{k+1} - \phi_k \right)  + h^{\gamma}\Gamma(2-\gamma) \sigma(W_{j-1}\phi_{j-1} + {\bf b}_{j-1}),\\ &\hspace{8cm}   2\leq j\leq L-1 , \\
				\phi_L   &= W_{L-1} \phi_{L-1} . 
			\end{aligned}
			\end{align}
			
			\item Adjoint Equation.
			
			\begin{equation}\label{eq:oc_adjoint}
				\begin{aligned}
					\psi_j &= \psi_{j+1} + \sum_{k=j+1}^{L-2} a_{k-j}\: \left(\psi_{k+1} - \psi_k \right)- \hspace{1cm} j = L-2,...,1\\
						&\hspace{2.5cm} h^{\gamma}\Gamma(2-\gamma)\left[-W_j^{T}\left(\psi_{j+1} \odot 	\sigma^{\prime}\left(W_{j} \phi_{j+1} + {\bf b}_j \right)\right)\right],\\
					\psi_{L-1}&=-W_{L-1}^T \psi_L, \\
					\psi_L &= \partial_{\phi_L}\mathcal{J(\bth ; \bxi,\bu)}.
			\end{aligned}
			\end{equation}				
			\item Derivative with respect to $\bth$.
			\begin{equation}\label{eq:design}
				\begin{aligned}
					\partial_{W_{L-1}}\mathcal{L}& =\psi_L\;\phi_{L-1}^T=\partial_{\phi_L}\mathcal{J(\bth ; \bxi,\bu)}\;\phi_{L-1}^T,\\
					\partial_{W_{j}}\mathcal{L}= &-\phi_{j}\:\left(\psi_{j+1} \odot \sigma^{\prime}(W_{j} \phi_{j} + {\bf b}_j)\right)^{T}+\partial_{W_j}\mathcal{J(\bth ; \bxi,\bu)},  \\
					&\hspace{5.5cm} j = 0,...,L-2,  \\
					\partial_{{\bf b}_{j}}\mathcal{L}= &-\; \psi_{j+1}^T\;\sigma^{\prime}(W_{j} \phi_{j} + {\bf b}_{j})+\partial_{{\bf b}_j}\mathcal{J(\bth ; \bxi,\bu)}.  \\
					&\hspace{5.5cm} j = 0,...,L-2 \, .
				\end{aligned}
			\end{equation}
		\end{enumerate}			
	\end{subequations}
The right-hand-side of \eqref{eq:design} represents the gradient of $\mathcal{L}$ with respect to $\bth$.	
We then use a gradient-based method (BFGS in our case) to identify $\bth$.

\subsection{Error analysis}
\label{erroranaly} 
In this section we briefly address the question of how well the deep neural approximation approximates the PDE solution map $\Phi:\mathbb{R}^{N _{\bxi}} \rightarrow \mathbb{R}^{N_x}.$ The approximation capabilities of neural networks has received a lot of attention recently in the literature, see e.g.,  \cite{SE94, Hornik94,  SX20, Yarotsky17} and the references therein. 
 The papers \cite{Hornik94, SX20} obtain results based on general activation functions
for which the approximation rate is $\mathcal{O}(n^{-1/2})$, where $n$ is the total number of hidden neurons. This implies  that neural networks can overcome the curse of dimensionality, as the approximation rate is independent of dimension $N_x$. 

We begin by making the observation that the fractional DNN can be written as a linear combination of activation functions evaluated at different layers. Indeed, observe from \eqref{eq:fRNN} that if $\Phi(\bxi)$  is approximated  by a one-hidden layer network $\hat\Phi(\bxi,\bth):=\phi_2$ (that is, $L=2$) then it  can be expressed as: 
\be
\label{2layer}
\phi_2= W_1\sigma (W_0\bxi +   {\bf b}_0).
\e
By a one-hidden layer network, we mean a network with the input layer, one hidden layer and the output layer \cite{SE94, Hornik94,  LLPS93, SX20}.
Next, 
 for  $L=3$, we obtain
\be
\label{phi3}
 \phi_3 = W_2\phi_2 = W_2\left[   \alpha_0\sigma(W_0\phi_0 + {\bf b}_0) + \alpha_1\sigma(W_1\phi_1 + {\bf b}_1) + \alpha_2\phi_0 \right],
\e
where $\alpha_0=1-a_1, \; \alpha_1=\tau:=h^{\gamma}\Gamma(2-\gamma)$  and $\alpha_2=a_1.$  Similarly, if we set $L=4$, then we get
\begin{eqnarray}
\label{phi4}
 \phi_4 &=& W_3\phi_3 \nonumber\\
   &=& W_3\left[   \alpha_0\sigma(W_0\phi_0 + {\bf b}_0) + \alpha_1\sigma(W_1\phi_1 + {\bf b}_1) + \alpha_2\sigma(W_2\phi_2 + {\bf b}_2) +\alpha_3\phi_0 \right],
\end{eqnarray}
where $\alpha_0=(1-a_1+ a_1^2 - a_2),\; \alpha_1= (1-a_1)\tau, \; \alpha_2=\tau,$  and $\alpha_3= 2a_2 - a_1 -a_1^2.$ 

Proceeding in a similar fashion yields the following result regarding multilayer fractional DNN.
\begin{proposition}[Representation of multilayer fractional DNN]
\label{multifdnn} 
For $L \ge 3$, the fractional DNN given by \eqref{eq:fRNN} fulfills
 \be
\label{phil}
\phi_L=W_{L-1}\left[\alpha_{L-1}\phi_0 + \sum_{i=0}^{L-2}\alpha_i\sigma(W_i\phi_i + {\bf b}_i )\right],
\e
where $\alpha_i$ are constants depending on $\tau$ and $a_i$, as given by \eqref{a_k}.
Observe that a one-hidden layer fDNN (that is, $L=2$) coincides with a one-hidden layer standard DNN. 
\end{proposition}

Error analysis for multilayer networks is generally challenging. Some papers that study this  include \cite{HLXZ18, Yarotsky17} and the references therein. The results  of these two papers focus mainly on the ReLU activation function. In particular, \cite{HLXZ18} discusses the approximation of linear finite elements by ReLU deep and shallow networks. 

In this paper, we restrict the analysis discussion to one-hidden layer fDNN i.e., $L=2$. 
In particular, we consider a one-hidden layer network with finitely many  neurons.
Notice that for $L=2$, fDNN coincides with standard DNN according to Proposition~\ref{multifdnn}. 
Therefore, it is possible to extend the result from \cite{SX20} to our case.  
To state the result from \cite{SX20},  we first introduce some notation. 

To this end, we make the assumption that  the function $\Phi(\bxi)$ is defined on a bounded domain $\mathcal{P}\subset \mathbb{R}^{N _{\bxi}}$ and has a bounded Barron norm 
\be
||\Phi||_{\mathcal{B}^m}=\int_{\mathcal{P}}(1 + \omega)^m |\hat{\Phi}(\omega)|\; d\omega, \;\;m\geq 0.
\e

Now, let $m\in \mathbb{N}\cup \{0\}$ and $p\in [1,\infty].$ Recall that the Sobolev space $\mathcal{W}^{m,p}(\mathcal{P})$
is the  space of functions  in $L^p(\mathcal{P})$ whose  weak derivatives of order $m$ are also in $L^p(\mathcal{P}):$ 
\be
\mathcal{W}^{m,p}(\mathcal{P}):=\{ f\in L^p(\mathcal{P}): D^{\bf m}f \in L^p(\mathcal{P}), \; \forall {\bf m} \;\mbox{with} \;|{\bf m}|\leq m\},
\e
where ${\bf m}=(m_1,\ldots,m_{N_{\xi}})\in \{0,1,\ldots \}^{N_{\xi}},$
$|{\bf m}|=m_1+m_2\ldots +m_{N_{\xi}}$  and $D^{\bf m}f$ is the weak derivative. In particular, the space $\mathcal{H}^m(\mathcal{P}):= \mathcal{W}^{m,2}(\mathcal{P})$
 is a Hilbert space with inner product and norm given respectively by
 \[
  (f,g)_{\mathcal{H}^m(\mathcal{P})}=\sum_{|{\bf m}|\leq m}\int_{\mathcal{P}} D^{\bf m}f({\bf x}) D^{\bf m}g({\bf x})\; d{\bf x},
 \]
and 
$ ||f||_{\mathcal{H}^m(\mathcal{P})} =(f,f)^{1/2}_{\mathcal{H}^m(\mathcal{P})}.$ 
For $p=\infty,$ we note that the Sobolev space $\mathcal{W}^{m,\infty}(\mathcal{P})$ is  equipped with the norm
\be
\label{infnorm}
 ||f||_{ \mathcal{W}^{m,\infty}(\mathcal{P})} =\max_{ {\bf m}:|{\bf m}|\leq m} \mbox{ess}\sup_{ {\bf x}\in \mathcal{P}}|D^{\bf m}f({\bf x})|.
\e
Also, we say that $f\in \mathcal{W}_{loc}^{m,p}(\mathcal{P})$ if $f\in \mathcal{W}^{m,p}(\mathcal{P}'),\; \; \forall $ compact ${\mathcal{P}'}\subset \mathcal{P}.$ 

Next, since the mapping $ \Phi:\mathbb{R}^{N_{\bxi}}  \rightarrow \mathbb{R}^{N_x}$ can be computed using ${N_x}$ mappings
$ \Phi_j:\mathbb{R}^{N_{\bxi}}  \rightarrow \mathbb{R},$ it suffices to focus on networks with one output unit. Observe from \eqref{2layer} that 
\be
\label{neurons}
\hat\Phi(\bxi,\bth):=\phi_2= \sum_{i=1}^n c_i \cdot \sigma (\omega_i\cdot \xi +   b_i),  
\e
where $n$ is the total number of neurons in the hidden layer, $b_i$ an element of the bias vector,  $c_i $ and $\omega_i$  are rows of $W_1$ and $W_0$ respectively. 
Now, for a given activation function $\sigma,$ define the set
\be
\label{Fn}
\mathcal{F}_{N_{\bxi}}^n(\sigma)=\left\{ \sum_{i=1}^n c_i\cdot\sigma (\omega_i\cdot\xi +   b_i), \quad  \omega_i\in\mathbb{R}^{N_{\xi}},\; b_i\in\mathbb{R} \right\}.
\e

We can now state the following result for any non-zero  activation functions satisfying some regularity conditions  \cite[Corollary 1]{SX20}. After this result, we will show that our activation function in \eqref{smrelu} fulfills the required assumptions.
\begin{theorem}
 Let  $\mathcal{P} $   be a bounded domain and  $\sigma \in \mathcal{W}_{loc}^{m,\infty}(\mathbb{R})$ be  a non-zero  activation function.
 Suppose there exists a function $\nu\in \mathcal{F}_1^{q}(\sigma)$ satisfying 
 \be
 \label{decay}
  |\nu^{(k)}(t)|\leq C_p(1 +  |t|)^{-p}, \;\; 0\leq k \leq m, \;\; p>1.
 \e
 Then, for any $\Phi\in \mathcal{B}^m,$ we have
 \be
 \label{errb}
 \inf_{\hat\Phi_n \in \mathcal{F}_{N_{\xi}}^n(\sigma) }|| \Phi - \hat\Phi_n||_{\mathcal{H}^m(\mathcal{P})} \leq C(\sigma,  m, p,\beta)|\mathcal{P}|^{1/2} {q}^{1/2}n^{-1/2}||\Phi||_{\mathcal{B}^{m+1}},
 \e
 where $\beta = \mbox{diam}(\mathcal{P})$ and $C$ is a constant depending on $\sigma , m, p$ and $\beta.$
\end{theorem}

Note that the bound in \eqref{errb} is $\mathcal{O}(n^{-1/2}).$
Here, the function $\nu(x)$ is a linear combination of the shifts and dilations of the activation $\sigma(x).$  Moreover, the activation function itself need not satisfy the decay  \eqref{decay}. The theorem says that it suffices to find some function $\nu\in \mathcal{F}_1^{q}(\sigma)$ (cf. (\ref{Fn}) with $n=q,\; N_{\xi}=1$) for which \eqref{decay} holds. 

\begin{figure}[h!]
		\centering
		\includegraphics[width=0.85\textwidth]{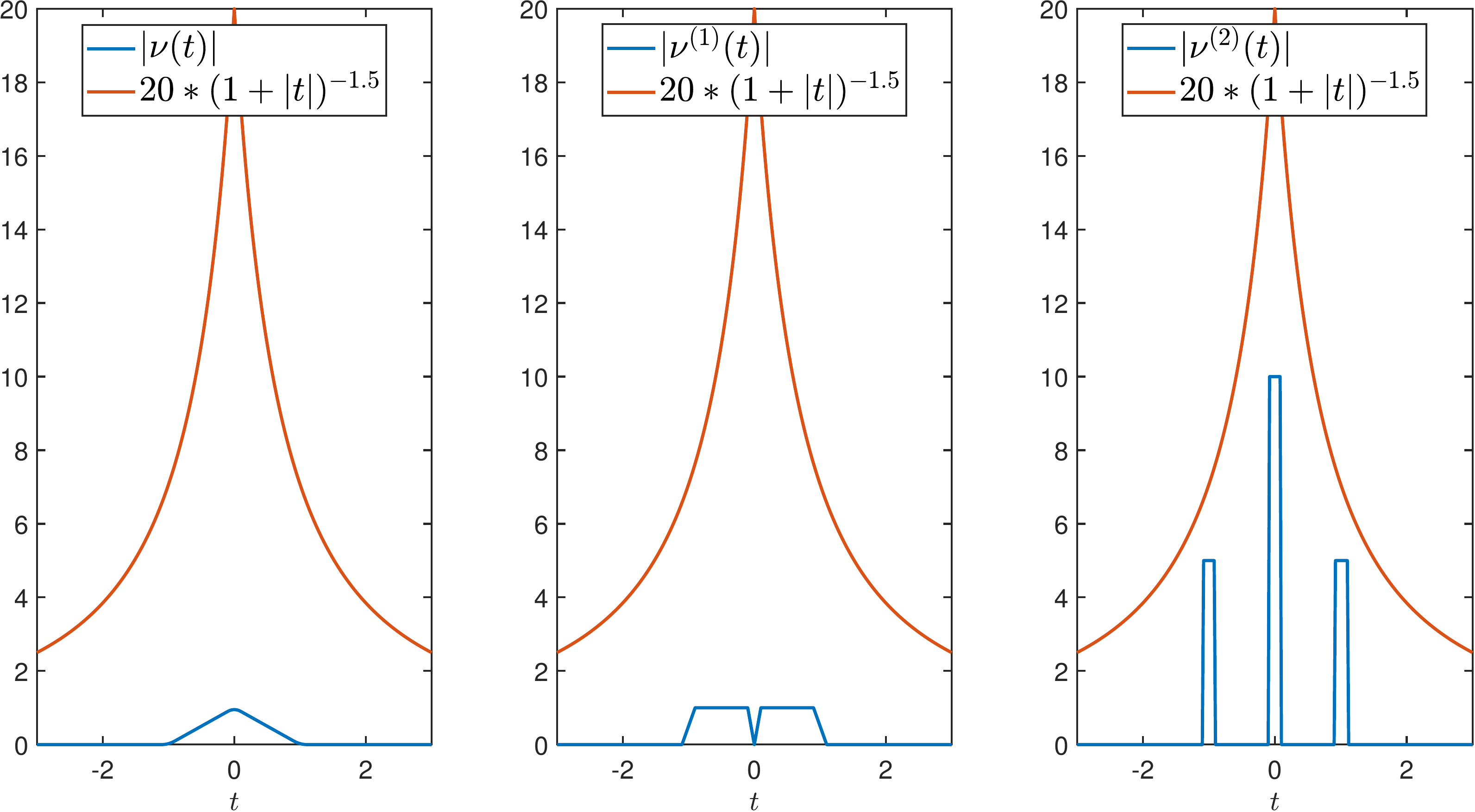}
		\caption{\label{fig:nu} Verification of the decay condition \eqref{decay} for $\nu$ in \eqref{eq:nu0}, which has been defined using the smooth ReLU in \eqref{smrelu}. }
	\end{figure}

In this paper, we are working with smooth ReLU (see \eqref{smrelu}) as the activation function $\sigma(x)$ for which $m=2$. One choice of $\nu$ for which \eqref{decay} holds is $\nu\in \mathcal{F}_1^3(\sigma)$ given by
	\begin{equation} 
		\label{eq:nu0}
		\nu(t)=\sigma(t+1) +\sigma(t-1)-2\sigma(t) .
	\end{equation}
	It can be shown that, $\nu$ satisfies \eqref{decay} with $C_p=20$ and $p=1.5$, that is,
	\be \label{eq:nu}
	 |\nu^{(k)}(t)|\leq 20(1 +  |t|)^{-1.5}, \;\; 0\leq k \leq 2. 
	\e This can be verified from Figure~\ref{fig:nu}. 
	
In what follows, we discuss the application of our proposed fractional DNN to the solution of two statistical Bayesian inverse problems.

\section{Application to Bayesian inverse problems}
\label{inverse_problems}
The  inverse
problem associated with  the forward problem (\ref{disc1}) essentially 
requires estimating a parameter vector $\bxi \in\mathbb{R}^{N_{\bxi}}$ given some observed noisy limited data
$d\in\mathbb{R}^{N_x}.$ 
In the Bayesian inference framework, 
the  posterior probability density,
$\pi(\bxi|d):\mathbb{R}^{N_{\bxi}}\rightarrow \mathbb{R},$   solves the
statistical inverse problem. In principle, 
$\pi(\bxi|d)$ encodes the uncertainty  from the set of observed data and the sought parameter vector.
More formally,  the posterior probability density  is given by  Baye's rule as
\begin{equation}
\label{post1}
\pi_{pos}(\bxi):=\pi(\bxi|d) 
     \propto  \pi(\bxi)  \pi(d|\bxi),
\end{equation}
where $\pi(\cdot):\mathbb{R}^{N_{\bxi}}\rightarrow \mathbb{R}$ is the prior and $\pi(d|\bxi)$ is the likelihood function. 
The standard approach to Bayesian inference uses the assumption that the observed data  are of the form (see e.g., \cite{CS07, Flath11})
\be
\label{obs}
 d=\Phi(\bxi) + \eta,
\e
where 
$\Phi$ is the  {\it parameter-to-observable} map
and $\eta \sim \mathcal{N}(0,\Sigma_{\eta})$. In our numerical experiments, we assume that $\Sigma_{\eta}=\kappa^2 I,$ where $I$ is the identity matrix with appropriate dimensions and $\kappa$ denotes the variance.
 Now, let the log-likelihood be given by
 \be
 \label{like}
 L(\bxi) = \frac{1}{2}\left|\left|\Sigma_{\eta}^{1/2}( d -\Phi(\bxi) )\right|\right|^2;
 \e
 then Baye's rule in (\ref{post1}) becomes
 \begin{equation}
\label{post2}
\pi_{pos}(\bxi):=\pi(\bxi|d) \propto  \pi(\bxi)  \pi(d|\bxi)
       =\frac{  1 }{ Z } \exp\left(-L(\bxi) \right)\pi(\bxi), 
     \end{equation}
where $Z=\int_{\mathcal{P}} \exp(-L(\bxi)) \pi(\bxi)\; d\bxi $  is a normalizing constant.

We note here that, in practice, the posterior density rarely admits analytic expression. Indeed, it is generally evaluated using Markov chain Monte Carlo (MCMC) sampling techniques.  In  MCMC schemes  one has to perform several thousands of
forward model simulations. This is a significant computational challenge, especially if the forward model  represents 
large-scale  high fidelity discretization of a nonlinear PDE. It is therefore reasonable to construct a cheap-to-evaluate surrogate for the forward model in the offline stage for use in online computations in MCMC algorithms. 

In what follows, we discuss   MCMC schemes for sampling   posterior probability distributions, as well as how our fractional DNN surrogate modeling strategy can be incorporated in them.

\subsection{Markov Chain Monte Carlo}
\label{dnn_mcmc} 
Markov chain Monte Carlo (MCMC) methods are very powerful and flexible numerical sampling approaches for approximating posterior distributions, see e.g., \cite{BGLFS17, CS07, RCS14}. Prominent MCMC methods include Metropolis-Hastings algorithm (MH)  \cite{RCS14}, Gibbs algorithm \cite{BSHL2014} and Hamiltonian Monte Carlo algorithm (HMC) \cite{BGLFS17, Bui14}, Metropolis-adjusted Langevin algorithm (MALA) \cite{BGLFS17}, and preconditioned Crank Nicolson algorithm (pCN) \cite{BSHL2014},  and their variants.
In our numerical experiments, we will use   MH,  pCN,  HMC, and MALA algorithms. We will describe only  MH in this paper and refer  the reader to  \cite{BGLFS17} for the details of the derivation and properties of the many variants of the other three algorithms.  

To this end, 
consider a given   parameter sample  $\bxi=\bxi^i$. The MH algorithm generates a new sample $\bxi^{i+1}$ as follows:
\begin{enumerate}
 \item[(1)] Generate a proposed sample $\bxi^{\ast}$ from a proposal density 
 $q(\bxi^{\ast}|\bxi )$, and compute $q(\bxi|\bxi^{\ast} )$.
 \item[(2)] Compute the acceptance probability
 \be
 \label{accpr}
 \alpha(\bxi^{\ast}|\bxi) =
 \min\left\{1, \frac{\pi(\bxi^{\ast}|d)q(\bxi^{\ast}|\bxi)}{\pi(\bxi |d)q(\bxi|\bxi^{\ast} ) } \right\}.
 \e
\item[(3)] If $\mbox{Uniform}(0; 1] <\alpha(\bxi^{\ast}|\bxi),$ then  $\bxi^{i+1} = \bxi^{\ast}.$
Else, set $\bxi^{i+1} = \bxi$.
\end{enumerate}
 Observe from  (\ref{like}) and (\ref{post2})  that  each evaluation of the likelihood function, and hence, the acceptance probability (\ref{accpr}) requires the evaluation of the forward model to compute $\pi(\bxi^{\ast}|d)$. In practice, one has to do  tens of thousands of forward  solves for the HM algorithm to converge. We propose to replace the forward solves by the fractional DNN surrogate model.  In our numerical experiments, we follow \cite{SHT01} in which an adaptive Gaussian proposal is used:
\be
\label{Gaussianp}
 q(\bxi^{\ast}|\bxi^i )=\exp\left(  -\frac{1}{2}(\bxi^{\ast} - \bxi^{i-1} )^TC_{i-1} (\bxi^{\ast} - \bxi^{i-1})\right),
\e
where $C_0=I,$ 
\[
 C_{i-1}=\frac{1}{i}\sum_{j=0}^{i-1}(\bxi^{j} -\bar\bxi )(\bxi^{j} -\bar\bxi)^T + \vartheta I,
\]
 $\bar\bxi = i^{-1}\sum_{j=0}^{i-1}\bxi^j,$ and $\vartheta \approx 10^{-8}.$
When the proposal is chosen adaptively as specified above, the HM method is referred to as an Adaptive Metropolis (AM) method \cite{JMB2018, SHT01}.

\section{Numerical Experiments}
\label{Numex}
In this section, we consider two statistical inverse problems. The first one is a diffusion-reaction problem in which two parameters need to be inferred \cite{EO19}. The second one is a more challenging problem -- a thermal fin problem from \cite{BGLFS17}, which involves one hundred parameters to be identified. All experiments were performed using MATLAB R2020b on a Mac desktop with RAM 32GB and CPU 3.6 GHz Quad-Core Intel Core i7. 

In both of these experiments we train using a 3-hidden layer network  with $15$ neurons in each hidden layer (that is, $L=4, \; n=45$) and $k=400$ neurons in the output layer, where $k$ matches  the dimension of POD basis as described in section \ref{subsec4a} below. We chose $\varepsilon=0.1$ in the smooth ReLU activation function, final time $T=1$ and step-size $h=\frac{1}{3}$. We set the fractional exponent $\gamma=0.5$ in the Caputo Fractional derivative and the regularization parameter $\lambda=10^{-6}$. To train the network, we use the BFGS optimization method \cite[Chapter 4]{CTK99}.

Tables \ref{train_bfgs_drp}  and \ref{train_bfgs_tfp} show the number of BFGS iterations and the CPU times required to train the data from the diffusion-reaction problem and the thermal fin problem, respectively. Also shown in these tables are the  relative errors obtained by evaluating the trained network at a parameter $\bxi^e$ not in the training set $\mathbb{E}$; here,
\[
 Error = \frac{||u(\bxi^e) -\hat\Phi(\bxi^e)||_{\infty}}{||u(\bxi^e)||_{\infty}},
\]
where $u(\bxi^e)$  and $\hat\Phi(\bxi^e)$ are the true and surrogate solutions at $\bxi^e$, respectively.
Note from both tables that after $1600$ BFGS iterations, the decrease in errors is not significant. Hence, for all the experiments discussed below,  we  use $1600$ BFGS iterations.

 \begin{table}[h!]
 \centering
\begin{tabular}{l|llll}
\hline
BFGS iterations    &  Error        &  Time  (in sec)        &  \\
\hline
\hline
$400$  &  $2.26\times 10^{-2}$  &    $ 7.99$    &  \\
\hline
$800$  &  $1.10\times 10^{-2}$  &    $ 15.10$    &  \\
\hline
$1600$ &  $4.09\times 10^{-3}$  &    $ 29.45$    & \\
\hline
$3200$ &  $3.43\times 10^{-3}$  &    $ 56.14$    & \\
\hline
$6400$ &  $2.56\times 10^{-3}$  &    $ 105.66$    & \\
\hline
\end{tabular}
\captionof{table}{Diffusion-reaction problem: Number of BFGS iterations, relative errors and times for training the fractional DNN.
}
\label{train_bfgs_drp}
\end{table}

 \begin{table}[h!]
 \centering
\begin{tabular}{l|llll}
\hline
BFGS iterations    &  Error        &  Time    (in sec)        &  \\
\hline
\hline
$400$  &  $1.87\times 10^{-2}$  &    $ 10.69$    &  \\
\hline
$800$  &  $1.15\times 10^{-2}$  &    $ 18.51$    &  \\
\hline
$1600$ &  $7.39\times 10^{-3}$  &    $ 33.52$    & \\
\hline
$3200$ &  $4.80\times 10^{-3}$  &    $ 62.90$    & \\
\hline
$6400$ &  $3.73\times 10^{-3}$  &    $ 117.55$    & \\
\hline
\end{tabular}
\captionof{table}{Thermal fin problem: Number of BFGS iterations, relative errors and times for training the fractional DNN.
}
\label{train_bfgs_tfp}
\end{table}

\subsection{Diffusion-reaction example}
\label{subsec4a} 
We consider the following nonlinear diffusion-reaction problem posed in a two-dimensional spatial domain \cite{Cha2011, GMNP07}
\begin{align}\label{nonlinear}
\begin{aligned}
-\Delta u + g(u;\bxi) &= f, \quad \mbox{in } \Omega = (0,1)^2, \\
		            u &= 0 , \quad \mbox{on } \partial\Omega , 
\end{aligned}		
\end{align}
where $g(u; \bxi) = \frac{\bxi_2}{\bxi_1}\left[\exp(\bxi_1 u) - 1\right]$ and $f = 100 \sin(2\pi x_1) \sin(2\pi x_2)$. 
Moreover, the parameters are $\bxi = (\bxi_1, \bxi_2)\in [0.01, 10]^2 \subset \mathbb{R}^2$.

Equation   (\ref{nonlinear}) is  discretized on a uniform mesh in $\Omega$ with  $64$ 
grid points in each direction using  centered differences resulting
in   $4096$ degrees of freedom. We obtained the  solution of the resulting 
system of nonlinear equations  using an inexact Newton-GMRES method 
as described in \cite{CTK95}. The  stopping tolerance was $10^{-6}.$

 To train the network, we first computed $N_s=900$ solution snapshots $\mathbb{S}$ corresponding to a set $\mathbb{E}$ of $900$ parameters $\bxi=(\bxi_1,\bxi_2)$ drawn from the parameter space $[0.01, 10]^2$. These parameters were chosen using  Latin hypercube sampling.
 Each solution snapshot is of dimension $N_x=4096.$
  Next,  we computed the SVD of the matrix of solution snapshots $\mathbb{S} = \widetilde{\mathbb{V}}\Sigma \mathbb{W}^T,$  and set our POD basis  $\mathbb{V}=\widetilde{\mathbb{V}}(:,1:k),$ where $k=400.$ Our training set then consisted of $\mathbb{E}$  as inputs and $\mathbb{V}^T\mathbb{S}\in \mathbb{R}^{k\times N_s}$  as our targets. As reported by Hesthaven and Ubbiali in \cite{HU2018} in the context of solving parameterized PDEs, the POD-DNN  approach accelerates the  online computations for surrogate models.
  Thus, we follow this approach in our numerical experiments; in particular, using  $k$-dimensional data (where $k\ll N_x$) confirms an overall speed up in the solution of the statistical inverse problems.

Next, in both numerical examples considered in this paper,  we solved  the inverse problems using $M=20,000$ MCMC samples.  In each case, the first $10,000$ samples were discarded for ``burn-in" effects, and the  remaining $10,000$ samples were used to obtain the reported statistics.   Here, we used an initial Gaussian 
proposal \eqref{Gaussianp} with covariance $C_0=I$ and updated the covariance after every $100$th sample. 
  As  in \cite{EO19}, we assume for this problem,
 that $\bxi$
is uniformly distributed over the parameter space $[0.01,10]^2$. Hence, 
 (\ref{post2})  becomes
\begin{align}
\label{postf3}
 \pi_{pos}(\bxi) &\propto
  \begin{cases}
     \exp\left(-\frac{1}{2\kappa^2}(d-\Phi(\bxi))^T (d-\Phi(\bxi)) \right)       & \text{if } \bxi \in [0.01,10]^2. \\
   0        & \text{otherwise.}
  \end{cases}
\end{align} 
We generated the observations $d$ by using 
 (\ref{obs}) the true parameter to be identified $\bxi^{e}=(1,0.1)$ and 
a Gaussian noise  vector $\eta$ with $\kappa =  10^{-2}.$

\begin{figure}[h!]
\centering
\includegraphics[width=0.8\textwidth,height=0.58\textwidth]{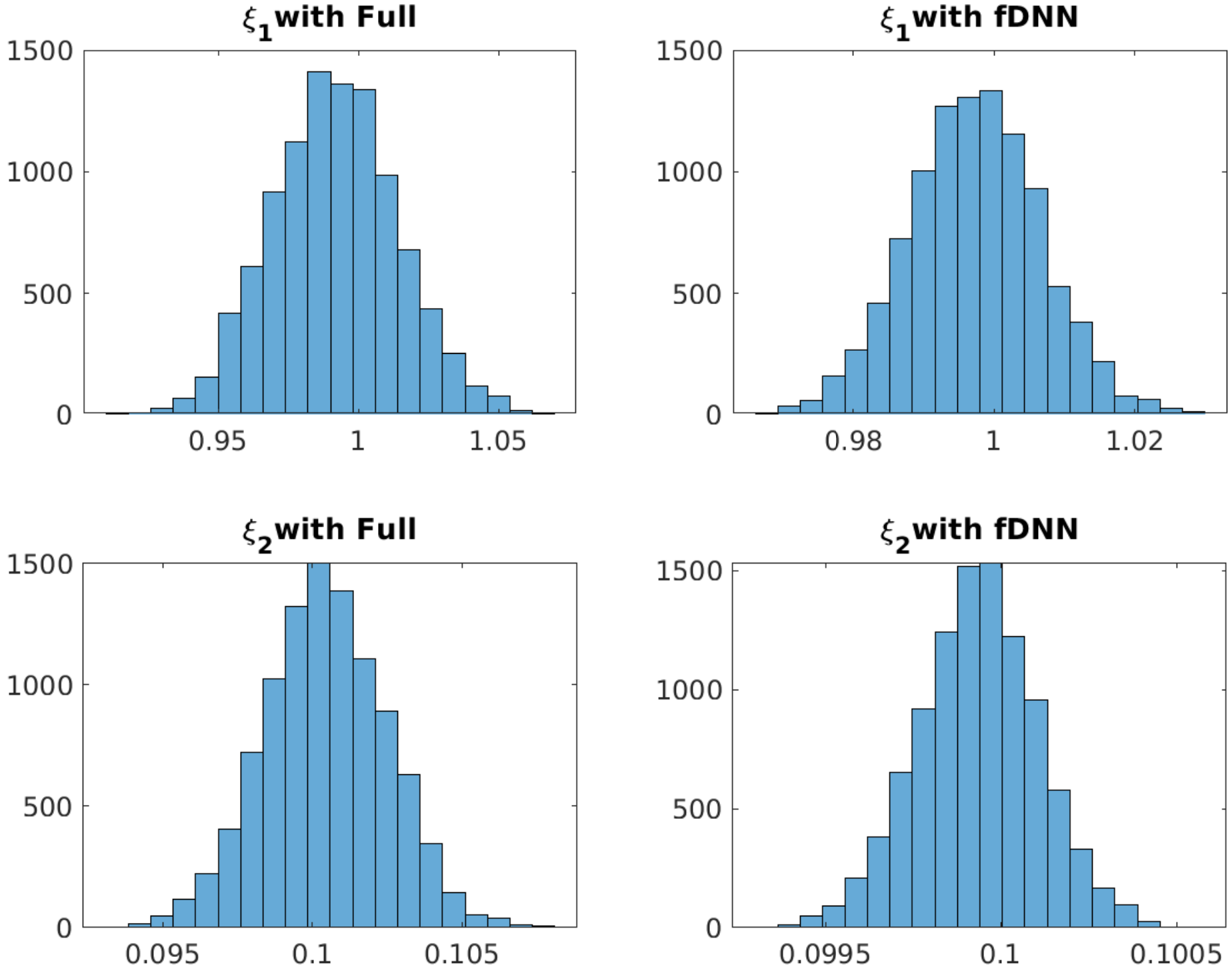}
\caption{Histograms  of  the  posterior  distributions  for the parameters $\bxi  = (\bxi_1, \bxi_2).$ 
They have been obtained from Full (left) and  fDNN (right) models with  $M=10,000$ MCMC samples.
}
\label{inv_hist101}
\end{figure}

\begin{figure}[h!]
\centering
\includegraphics[width=0.9\textwidth,height=0.75\textwidth]{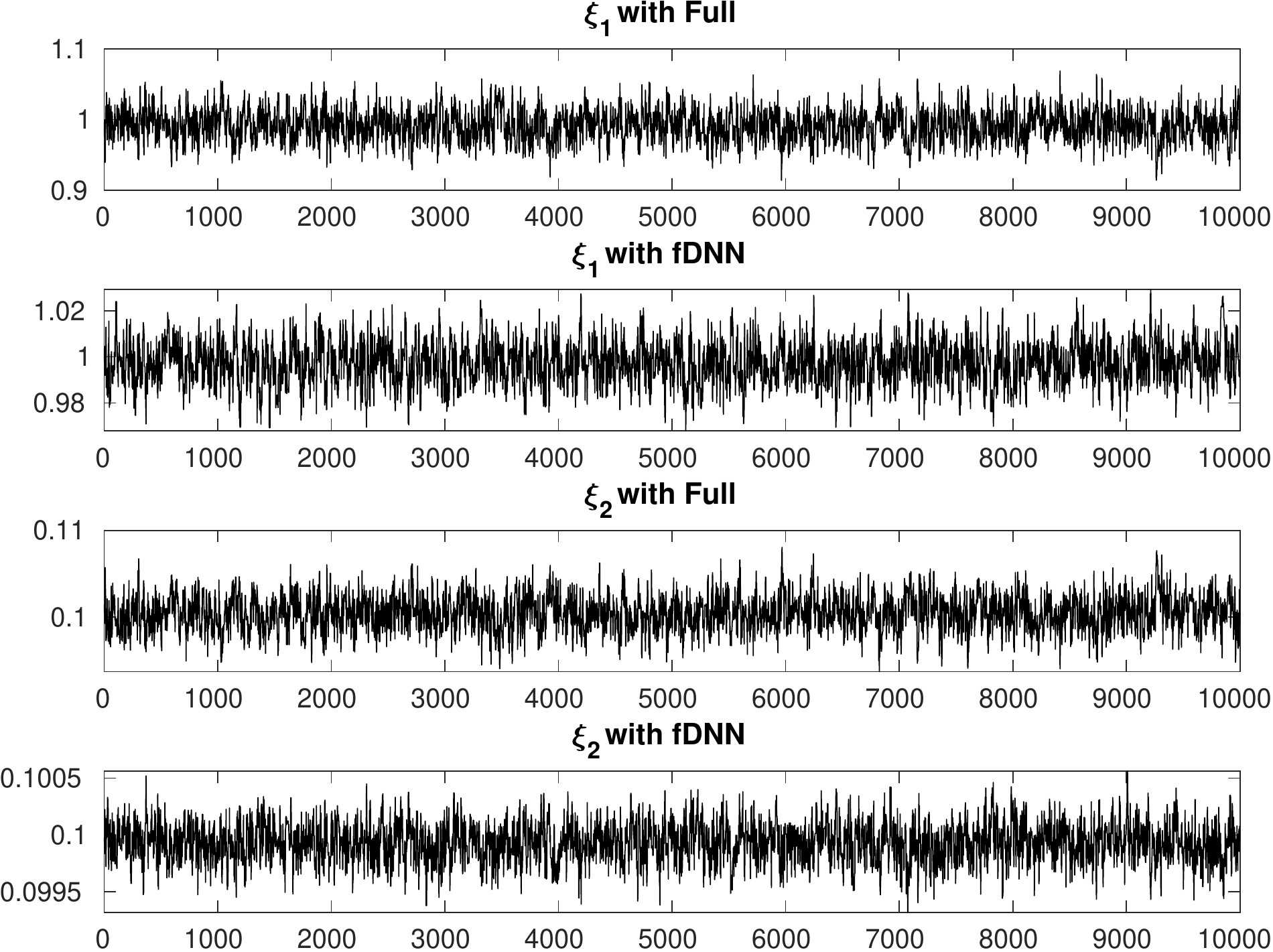}
\caption{MCMC samples
for the parameters $\bxi  = (\bxi_1, \bxi_2)$ using Full (first and third) and  fDNN (second and fourth) models. For the full model, the $95\%$ confidence interval for $\xi_1$ and $\xi_2$ are  $[0.9496, 1.0492]$ and $[0.0954, 0.1048]$. For the fDNN model, 
the $95\%$ confidence interval for $\xi_1$ and $\xi_2$ are 
$[0.9818, 1.0196]$ and $[0.0995, 0.1005]$.
}
\label{inv_mc101}
\end{figure}

\begin{figure}[h!]
\centering
\includegraphics[width=1.0\textwidth,height=0.4\textwidth]{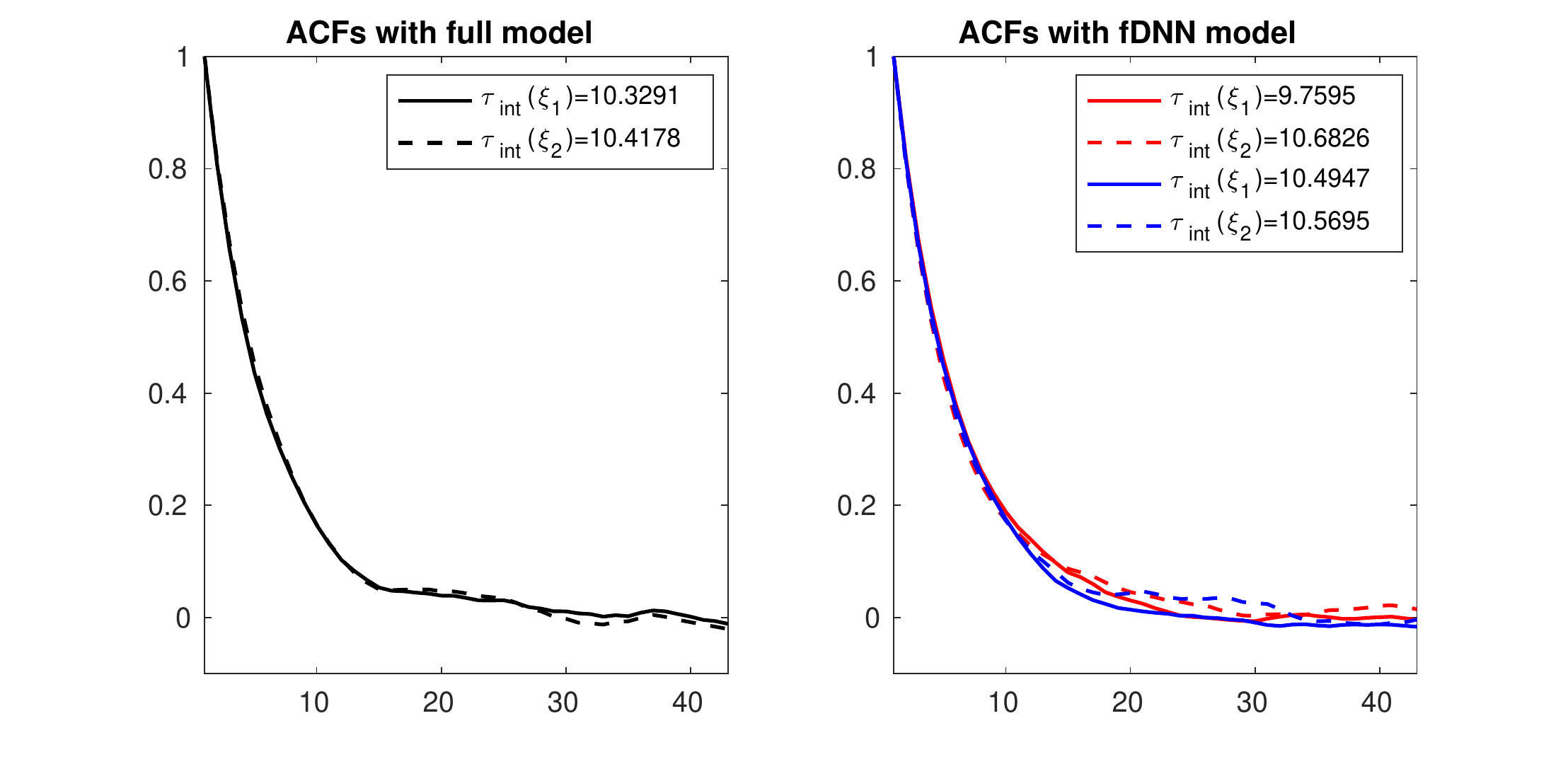}
\caption{Autocorrelation functions (ACFs)  for 
  $\bxi_1$ and $\bxi_2$ chains computed with Full (left) and  DNN (right) models. The red and blue lines are ACFs obtained with 1600 and 6400 BFGS iterations, respectively.}
\label{inv_acf101}
\end{figure}

   Figures \ref{inv_hist101},  \ref{inv_mc101}  and \ref{inv_acf101}, represent, respectively,
 the histogram (which depicts the posterior distribution), the Markov chains  and the autocorrelation functions corresponding to the parameters using the high fidelity model (Full) and the deep neural network surrogate (DNN) models in the MCMC algorithm. Recall that, for 
a Markov chain $\{\delta_j\}_{j=1}^J$ generated by the Metropolis-Hastings algorithm,
 with variance $\kappa^2,$
the autocorrelation function (ACF) $\varrho$ of the  $\delta$-chain is given by
\[
 \varrho(j) = \mbox{cov}(\delta_1,\delta_{1+|j|})/\kappa^2. 
\]
 and the integrated autocorrelation time, $\tau_{int}$, (IACT)
of the chain is defined as
\begin{eqnarray}
 \tau_{int}(\delta)&:=&\sum_{j=-J+1}^{J-1}\varrho(j)
 \approx 1 + 2\sum_{j=1}^{J-1}\left(1-\frac{j}{J}\right)\mbox{cov}(\delta_1,\delta_{1+j})/\kappa^2.
 \end{eqnarray}
ACF decays very fast to zero as $J\rightarrow \infty.$ In parctice, $J$ is often taken 
to be $\left\lfloor{10\log_{10}M}\right\rfloor,$ \cite{BSHL2014}.
If $\mbox{IACT} = K,$  this means that the roughly every $K$th sample 
of the $\delta$ chain is independent.  
We have used the following  estimators to approximate $\varrho(j)$  and $\tau_{int}$ in our computations \cite{JMB2018, Sokal97}: 
\[
  \varrho(j):=B(j)/B(0), \;\; \mbox{and}\;\; \tau_{int} :=  \sum_{j=-\widehat J}^{\widehat J}\varrho(j),
\]
where 
$B(j)=\frac{1}{J-j}\sum_{k=1}^{J-j} (\delta_k - \bar\delta)(\delta_{k+|j|}- \bar\delta),$
 $\bar\delta$ is the mean of the  $\delta$-chain, and ${\widehat J}$ is chosen be the smallest integer such that ${\widehat J} \geq 3\tau_{int}.$

Observe from Figs. \ref{inv_hist101} and \ref{inv_mc101} that, for both Full and fDNN models, the respective histograms  and Markov chains are centered around the parameters of interest $(1,0.1).$  In fact, for the full model, the $95\%$ confidence intervals (CIs) for $\xi_1$ and $\xi_2$ are  $[0.9496, 1.0492]$ and $[0.0954, 0.1048]$. For the fDNN model, the CIs  are $[0.9818, 1.0196]$ and $[0.0995, 0.1005]$. Thus,  fDNN identifies $(1,0.1)$ appropriately. 

We also examine the impact of training accuracy in  Fig \ref{inv_acf101}.  The figure  shows that the ACFs decay very fast to zero, as expected. The value of IACT is roughly $10$ for all the   MCMC chains generated by the full model (left). The right of the figure shows the IACT for fDNN model, where the red-colored curves show the results obtained with $1600$ training steps and the blue curves show those obtained with $6400$ steps. In general, the IACT values imply that roughly every  $10th$  sample of the MCMC chains  generated by both the full and the fDNN models are independent. 

These results indicate that the fDNN surrogate model produces results with accuracy 
comparable to those obtained using the full model.  The advantage of  the surrogate lies in
its reduced costs. 
In this example, using the HM algorithm, the  full model required $2347.7$ seconds of CPU
time to solve the inverse problem, whereas the fractional DNN model required $29.3$ 
seconds (the online costs for fDNN), a reduction in CPU time of a factor of $80$.
As always for surrogate approximations, there is overhead, the offline computations, 
associated with construction of the surrogate, i.e.,  
identification of the parameter set $\bth$ that defines the fDNN.
For this example, this entailed construction of the $900$ snapshot solutions used to
produce targets for the training process, computation of the SVD of the matrix $\mathbb{S}$ 
of snapshots giving the targets, together with the training of the network (using BFGS).
The times required for these computations were $127.8$ seconds to construct the inputs/targets,  $0.52$ second to compute the SVD and  $29.5$ seconds to train the network.
Thus, the offline and online times sum to $187.1$ seconds, which is a lot smaller than the
time  needed for the full solution  (that is, $2347.7$ seconds).
Of course, the offline computations represent a one-time expense which need not be repeated
if more samples are used to perform an MCMC simulation or if different MCMC algorithms 
are used.\footnote{We did not try to optimize the number of samples used for training and in particular $N_s = 900$  was an essentially arbitrary choice. The results obtained here with $N_s = 900$ were virtually the same as those found using $10,000$ training samples, with the smaller number of samples incurring dramatically lower offline costs.
Further reduction in computational time could be achieved using a smaller training data set.
In \cite{HU2018}, the authors used $\mathcal{O}(400)$ snapshots to train the feedforward network and still achieved good results.
} 
This is, for instance, the case with the Differential Evolution Adaptive Metropolis  (DREAM) method which runs multiple different chains simultaneously when used to  sample the posterior distributions \cite{Vrugt2009}.

\subsection{Thermal fin example}
\label{subsec4b} 
\begin{figure}[h!]
\centering
\includegraphics[width=1.0\textwidth]{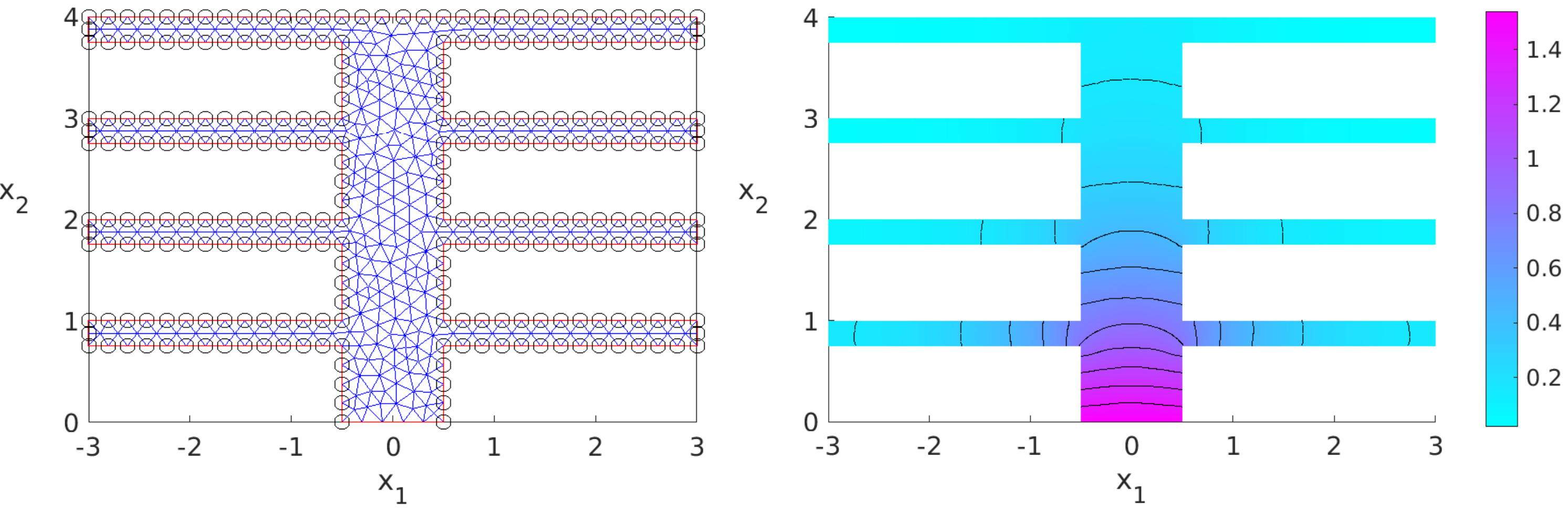}
\caption{The location of observations (circles) (left) and the forward PDE solution $u$ under the true parameter $\bxi$ (right).
}
\label{forwsol}
\end{figure}

 Next, we consider the following thermal fin problem from \cite{BGLFS17}:
\begin{align}
 \label{soc}
 \begin{aligned} 
    -\mbox{div } (e^{\bxi({\bf x})}\nabla u )  &=0,  \quad  &&\mbox{in } \Omega , \\
   (e^{\bxi({\bf x})}\nabla u )\cdot {\bf n} + 0.1 u &= 0, \quad &&\mbox{on }  \partial\Omega \setminus\Gamma,\,\\
   (e^{\bxi({\bf x})}\nabla u({\bf x}) )\cdot {\bf n} &= 1, \quad &&\mbox{on } \Gamma=(-0.5,  0.5) \times \{ 0 \}.
\end{aligned}
\end{align}
These equations (\ref{soc}) represent a forward model for  heat conduction over the non-convex domain $\Omega$ as shown in Fig. \ref{forwsol}.
Given the heat conductivity function $e^{\bxi({\bf x})}$, the forward problem (\ref{soc}) is used to compute the temperature $u$. The goal of this inverse problem is to  infer $100$ unknown parameters
$\bxi$ from  $262$  noisy observations of $u$. Fig. \ref{forwsol} shows the location of the  observations on the boundary $\partial\Omega \setminus\Gamma$, as well as  the forward PDE solution $u$ at the true parameter $\bxi$.

To train the network,   we  first computed, as in the diffusion-reaction case,  $N_s=900$ solution snapshots $\mathbb{S}$ corresponding to $900$ parameters $\bxi\in \mathbb{R}^{100}$ drawn  using  Latin hypercube sampling. This problem is a lot more difficult than the previous problem in the sense that the dimension of the parameter space in this case is $100$ rather than  $2.$ Next, we compute the SVD of $\mathbb{S}$ and proceed as before. 

\begin{figure}[h!]
\centering
\includegraphics[width=1.0\textwidth,height=0.65\textwidth]{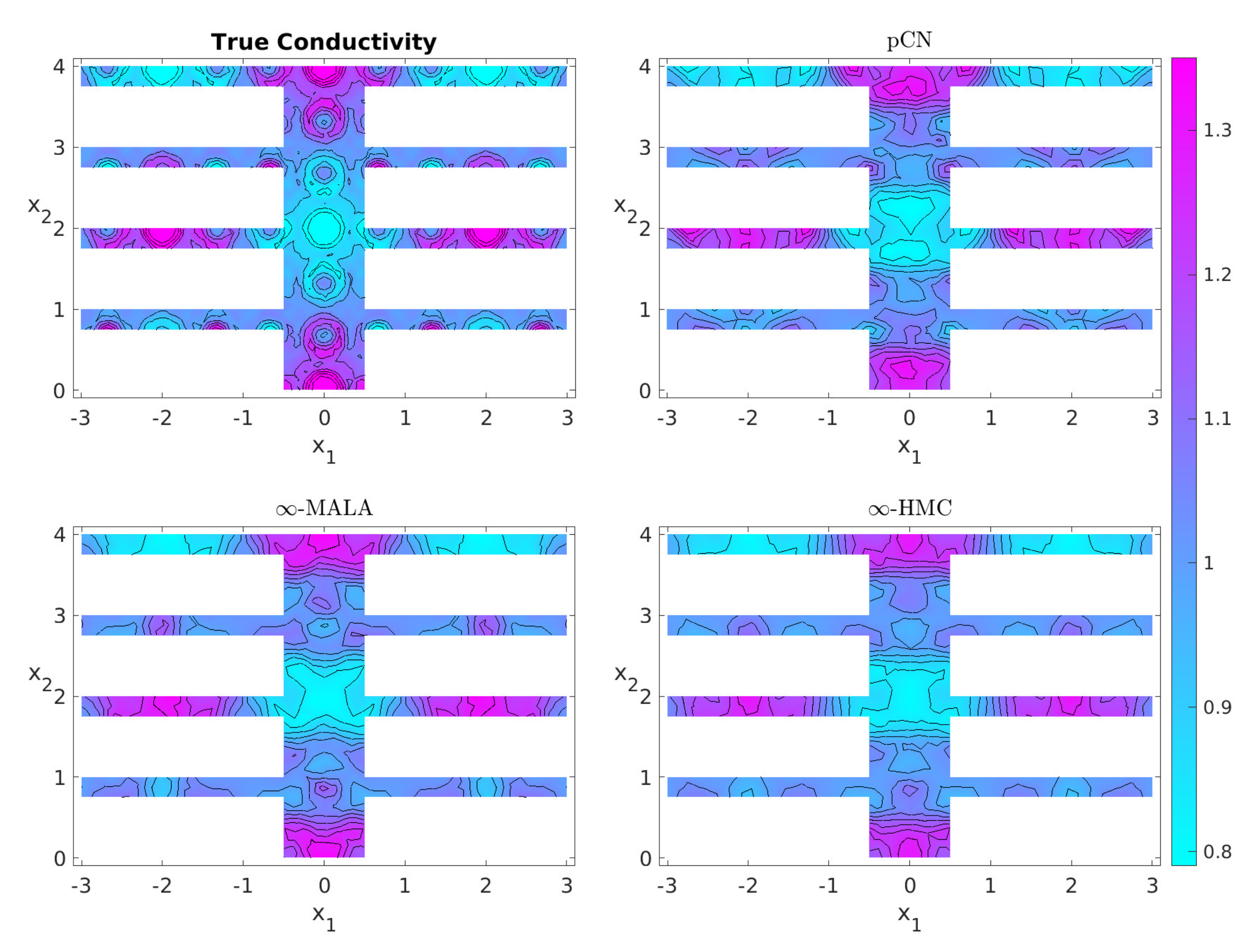}
\caption{The true heat conductivity field $e^{\bxi({\bf x})}$
(upper left) and the mean estimates of the posterior  obtained by the different MCMC methods using  fractional DNN as a surrogate model. }
\label{truth_mcmc}
\end{figure}

In what follows, the infinite variants of Riemannian manifold Metropolis-adjusted Langevin algorithm (MALA)  and  Hamiltonian Monte-Carlo (HMC) algorithms as presented in \cite{BGLFS17} are denoted, respectively, by $\infty$-MALA and  $\infty$-HMC. In particular,  Fig. \ref{truth_mcmc} shows the true conductivity, as well as the  posterior mean estimates obtained by three  MCMC methods -- pCN, $\infty$-MALA and  $\infty$-HMC -- using the fractional DNN as a surrogate model. 

As in \cite{BGLFS17}, for this problem we used a Gaussian prior  defined on the domain $\mathcal{D}:=[-3,3]\times [0,4]$ with covariance $\mathcal{C}$ of eigen-structure given by
\[
\mathcal{C} = \sum_{i\in I}\mu_i^2\{ \varphi_i \otimes \varphi_i\},
\]
where
$\mu_i^2 =\{ \pi^2( (i_1 + 1/2)^2 + (i_2 + 1/2)^2  )  \}^{-1.2},$
\[
\varphi_i({\bf x}) =  2|\mathcal{D}|^{-1/2}\cos(\pi(i_1 + 1/2)x_1)\cos(\pi(i_2 + 1/2)x_2),
\]
and
$I = \{i=(i_1,i_2), i_1\geq 0, \; i_2\geq 0\}. $
Moreover,  the log-conductivity $\bxi$ was the true parameter  to be inferred  with coordinates
$\bxi_i = \mu_i \sin (( i_1 - 1/2 )^2 + ( i_2 - 1/2 )^2 ), \;\; i_1 \leq 10,\; i_2 \leq 10.$

Table \ref{acc_time} shows the average acceptance rates for these models and the computational
times required to perform the MCMC simulation for different variants of the MCMC algorithm,
using both fDNN  surrogate computations and full-order discrete PDE solution.
These results indicate that  the acceptance rates are comparable for the fDNN surrogate model and the full order forward discrete PDE solver. 
Moreover, there is about a $90\%$ reduction in the computational times when fDNN surrogate
solution is used, a clear advantage of the fDNN approach.
The costs of the offline computations for fDNN were 
$63.2$ seconds to generate the data used for the fractional DNN models and  $33.5$ seconds to train the network (with $1600$ BFGS iterations) a total of $99.7$ seconds. 
As  in the case of  the diffusion-reaction problem, the overhead required for the 
one-time offline computations is offset by the short CPU time required for the online phase 
with the MCMC algorithms.
Note that the same offline computations were used for each of the three MCMC simulations tested.

 \begin{table}[h!]
 \centering
\begin{tabular}{l|lllll}
\hline
Model    &  pCN        &  $\infty$-MALA          &  $\infty$-HMC &  \\
\hline
\hline
Acc. Rate (fDNN)  &  $0.67$  &   $0.67$    & $0.79$  & \\
\hline
Acc. Rate (Full)  &  $0.66$  &   $0.70$    & $0.75$  & \\
\hline
\hline
CPU time  (fDNN) &  $16.46$  &    $ 98.0$    & $228.4$   & \\
\hline
CPU time  (Full) &  $157.8$  &    $ 958.9$    & $2585.3$   & \\
\hline
\end{tabular}
\captionof{table}{Acceptance rates and computational times  
needed 
 to solve the inverse problem by pCN, $\infty$-MALA and $\infty$-HMC algorithms together with fDNN  and  full forward  models.
}
\label{acc_time}
\end{table}

\section{Conclusions}
\label{conclusions} 
This paper has introduced a novel deep neural network (DNN) based approach to approximate 
nonlinear parametrized partial differential equations (PDEs). The proposed DNN helps learn the 
parameter to PDE solution operator. We have used this learnt solution operator to solve two challenging 
Bayesian inverse problems. The proposed approach is highly efficient without
compromising on  accuracy. We emphasize, that the proposed approach shows several 
advantages over the traditional surrogate approaches for parametrized PDEs 
such as reduced basis methods. 
For instance the proposed approach is fully non-intrusive and therefore it can be directly 
used in legacy codes, unlike the reduced basis method for nonlinear PDEs which can be 
highly intrusive.

\bibliographystyle{siam}
\bibliography{refs}
\end{document}